\documentclass[pdflatex,sn-mathphys-num]{sn-jnl}

\usepackage{amsmath,amssymb,amsthm,mathtools}
\usepackage{mathrsfs}

\theoremstyle{thmstyleone}
\newtheorem{theorem}{Theorem}
\newtheorem{lemma}[theorem]{Lemma}
\theoremstyle{thmstyletwo}

\theoremstyle{thmstylethree}

\newcommand{\co}{\operatorname{co}}
\newcommand{\argmin}{\operatorname*{arg\,min}}
\newcommand{\R}{\mathbb{R}}
\newcommand{\Q}{\mathbb{Q}}

\begin{document}

\title{A Fermat--Torricelli convex-hull characterization
	of inner-product spaces via Birkhoff--James orthogonality}

\author[1]{\fnm{Zhiyao} \sur{Fang}}
\email{fzy73@stu.aqnu.edu.cn}

\author[1]{\fnm{Junxiang} \sur{Qi}}
\email{y25060009@stu.aqnu.edu.cn}

\author*[1]{\fnm{Qi} \sur{Liu}}
\email{liuq67@aqnu.edu.cn}

\affil[1]{%
    \orgdiv{School of Mathematics and Statistics},
    \orgname{Anqing Normal University},
    \orgaddress{
        \city{Anqing},
        \postcode{246133},
        \state{Anhui},
        \country{China}
    }
}

\abstract{The classical location problem posed by Weber
	[Über den Standort der Industrien,
		Tübingen: J. C. B. Mohr, 1909]
	has motivated various geometric characterizations
	of inner-product spaces. Ben\'itez, Fern\'andez and Soriano
	[Trans. Amer. Math. Soc. \textbf{354} (2002), 5027--5038]
	characterized real inner-product spaces through the
	convex-hull property of Fermat--Torricelli medians.
	Pakhrou [J. Funct. Anal. \textbf{285} (2023),  110078]
	further established a characterization using $\gamma$-centers
	and equilateral triangles. Building on these results, we
	show that the Fermat--Torricelli convex-hull criterion
	remains valid when restricted to triples satisfying a
	one-sided Birkhoff--James orthogonality condition.
	The proof combines norming functionals, radial rescaling,
	Alexandrov differentiability, and the area formula.
	No completeness, smoothness, or strict convexity
	assumptions are required.}

	
\maketitle

\section{Introduction}\label{sec:introduction}

For points \(a_1,a_2,a_3\) in a normed space \(X\), their
Fermat--Torricelli median set\footnote{In a general normed space, the set of
	Fermat--Torricelli medians may be empty or contain
	more than one point; see the existence and uniqueness
	examples in \cite[pp.~5027--5028]{BFS2002}.}is
\[
Z_1(a_1,a_2,a_3)
=\argmin_{z\in X}\bigl(\|a_1-z\|+\|a_2-z\|+\|a_3-z\|\bigr).
\]
The convex-hull property asks whether this set meets \(\co\{a_1,a_2,a_3\}\). Ben\'itez, Fern\'andez and Soriano proved that, for a real normed space of dimension at least three, this property holds for every triple precisely when the norm is induced by an inner product \cite[Theorem~15]{BFS2002}. The hypothesis there tests all triangles.

\begin{theorem}
	\cite[Theorem~15]{BFS2002}
	Let \(X\) be a real normed space with \(\dim X\geq 3\).
	Then \(X\) is an inner-product space if and only if
	\[
	Z_1(u,v,w)\cap\operatorname{co}\{u,v,w\}
	\ne\varnothing,
	\qquad \text{for all }u,v,w\in X.
	\]
\end{theorem}

Related convex-hull characterizations were obtained by
Ben\'itez, Fern\'andez and Soriano for 2-centers and
weighted $p$-centers
\cite{Benitez2Centers2000}.
Ben\'itez and Y\'a\~nez subsequently established further
characterizations of inner-product spaces involving
middle points and Fermat--Torricelli medians
\cite{BenitezYanez2007}.
They also obtained a characterization based on the
midpoint properties of chords of the unit sphere
\cite{BenitezYanezBLMS2007}.
More recently, Pakhrou \cite{Pakhrou2023} extended
the convex-hull approach to $\gamma$-centers associated
with symmetric monotone norms, showing that it suffices
to consider equilateral triangles.
In the present paper, we strengthen the Fermat--Torricelli
convex-hull characterization by restricting the test to
triangles satisfying a single Birkhoff--James
orthogonality condition.

Pakhrou \cite[Theorem~5.4 and Corollaries~5.5--5.6]{Pakhrou2023}
characterized real inner-product spaces of dimension at least
three in terms of $\gamma$-centers, where $\gamma$ is a
symmetric monotone norm on $\mathbb{R}^3$, and showed that
it suffices to consider equilateral triangles. In particular,
taking $\gamma=\|\cdot\|_1$ yields the corresponding
Fermat--Torricelli center characterization. In contrast,
our characterization restricts the test to triples
$\{0,x,y\}$ satisfying the one-sided Birkhoff--James
orthogonality condition $x\perp_B y$, thus providing a
different orthogonality-based criterion for inner-product
spaces.While Pakhrou's characterization relies on the equality
of the three side lengths of an equilateral triangle,
our approach is based on the inequality defining
Birkhoff--James orthogonality. These two restrictions
reflect different geometric properties and lead to
distinct criteria for characterizing inner-product spaces.

Birkhoff--James orthogonality is a rich and interesting
notion with applications in various areas; see
\cite{BhatiaSemrl1999,
	ZamaniMoslehian2016,MartinEtAl2024,Khurana2024,
	ChengFelder2023}
and the references therein.

Our proof combines radial rescaling and norming functionals
(convex analysis; Lemma~\ref{lem:radial}) with a
Lipschitz-graph cover, Alexandrov's theorem, and the area
formula (geometric measure theory; Sections~\ref{sec:nonsmooth}
and~\ref{sec:selection}) to handle nonsmooth norms
.

\begin{theorem}\label{thm:main}
Let \(X\) be a real normed space with \(\dim X\geq 3\). The following are equivalent:
\begin{enumerate}
\item the norm of \(X\) is induced by an inner product;
\item for every \(x,y\in X\),
\[
x\perp_B y
\quad\Longrightarrow\quad
Z_1(0,x,y)\cap\co\{0,x,y\}\ne\varnothing.
\tag{1.1}\label{eq:restricted-property}
\]
\end{enumerate}

\end{theorem}

\section{Properties of Separated Triples}
\label{sec:rigidity}

We first work in dimension three. Fix an auxiliary Euclidean structure on \(\R^3\). Write \(|\cdot|\) for the Euclidean norm and \(N\) for the norm under consideration. All Lebesgue measure and differentiability statements below refer to these coordinates. For \(u\ne0\), let
\[
J(u)=\{f\in(\R^3)^*: \|f\|_{*}=1,\ f(u)=N(u)\}
\]
be the set of norming functionals. It is a nonempty
compact convex set, \(J(tu)=J(u)\) for \(t>0\), and
its graph is closed. By James' characterization
\cite{James1947},
\[
x\perp_B y
\quad\Longleftrightarrow\quad
\text{there exists }f\in J(x)\text{ with }f(y)=0
\qquad (x\ne0).
\tag{2.1}\label{eq:bj-support}
\]
For nonzero \(u_1,u_2,u_3\), the optimality condition for the Fermat-Torricelli problem is
\[
0\in Z_1(u_1,u_2,u_3)
\quad\Longleftrightarrow\quad
\text{there are }f_i\in J(u_i)\text{ such that }f_1+f_2+f_3=0.
\tag{2.2}\label{eq:ft-balance}
\]
This is also Proposition~1 of \cite{BFS2002}.

For a unit dual functional \(f\), define its contact cone by
\[
C_f=\{v\in\R^3:f(v)=N(v)\}.
\]
If \(f_i\in J(u_i)\) and \(\sum_i f_i=0\), then
\[
Z_1(u_1,u_2,u_3)=\bigcap_{i=1}^3(u_i-C_{f_i}).
\tag{2.3}\label{eq:median-contact}
\]
Indeed, the three nonnegative quantities \(N(u_i-z)-f_i(u_i-z)\) sum to
\[
\sum_iN(u_i-z)-\sum_iN(u_i),
\]
so their sum vanishes exactly when each vanishes. In particular, (2.3) remains valid when one of the residual vectors is zero.

We call a residual triple \(u=(u_1,u_2,u_3)\) separated if, for a fixed nonzero linear functional \(p\),
\[
0\in Z_1(u_1,u_2,u_3),\qquad
p(z)\geq0\ \ (z\in Z_1(u_1,u_2,u_3)),\qquad
p(u_i)<0\ \ (i=1,2,3).
\tag{2.4}\label{eq:S}
\]
Such a triple has its entire median set disjoint from the convex hull of its vertices.

\begin{lemma}\label{lem:radial}
If \(u\) satisfies \eqref{eq:S}, then \((a_1u_1,a_2u_2,a_3u_3)\) also satisfies \eqref{eq:S} for every \(a_1,a_2,a_3>0\).
\end{lemma}

\begin{proof}
	By Corollary~3 of \cite{BFS2002}, the origin remains a
	median of the rescaled triple. Choose \(f_i\in J(u_i)\)
	with \(\sum_{i=1}^3 f_i=0\).
	
	Let \(z\) be any median of the rescaled triple. By
	\eqref{eq:median-contact}, \(a_i u_i-z\in C_{f_i}\).
	Choose \(\delta>0\) such that \(\delta a_i\leq1\).
	Since \(C_{f_i}\) is a convex cone,
	\[
	u_i-\delta z
	=(1-\delta a_i)u_i+\delta(a_i u_i-z)
	\in C_{f_i}.
	\]
	Hence \(\delta z\in Z_1(u_1,u_2,u_3)\), so
	\(p(\delta z)\geq0\). Thus \(p(z)\geq0\).
	Since \(p(a_i u_i)<0\), the claim follows.
\end{proof}

\begin{lemma}\label{lem:mixed-edge}
Assume that \(N\) satisfies the restricted property \eqref{eq:restricted-property}. If \(u\) satisfies \eqref{eq:S}, then, for distinct \(i,j,k\) and all \(s,t>0\),
\[
J(su_j-tu_k)\subset\{\varphi:\varphi(u_i)=0\},
\tag{2.5}\label{eq:annihilation}
\]
and every \(\varphi\in J(su_j-tu_k)\) satisfies
\[
\varphi(u_j)>0,\qquad \varphi(u_k)<0.
\tag{2.6}\label{eq:strict-signs}
\]
\end{lemma}

\begin{proof}
Write the residuals as \(a,b,c\). For \(\alpha,\beta,\gamma>0\), define the compact interval
\[
I(\alpha,\beta,\gamma)
=\{\varphi(\gamma c-\alpha a):\varphi\in J(\beta b-\alpha a)\}.
\]
If this interval contained zero, then there would exist
\(\varphi\in J(\beta b-\alpha a)\) such that
\[
\varphi(\gamma c-\alpha a)=0.
\]
By \eqref{eq:bj-support}, it follows that
\[
\beta b-\alpha a\perp_B\gamma c-\alpha a.
\]The restricted property would give a median in the convex hull of \(\{\alpha a,\beta b,\gamma c\}\), contrary to Lemma~\ref{lem:radial}. Thus \(I\) avoids zero.

The map \(u\mapsto J(u)\) is upper semicontinuous away from the origin. It follows that the sets of parameters for which \(I\) is strictly positive or strictly negative are open. The positive octant is connected, and \(I\) never contains zero, so its sign is constant there. For fixed \(\beta,\gamma\), and all sufficiently large \(\alpha\), every \(\varphi\in J(\beta b-\alpha a)\) satisfies
\[
\varphi(\gamma c-\alpha a)
=\varphi(\beta b-\alpha a)+\varphi(\gamma c-\beta b)
\geq N(\beta b-\alpha a)-N(\gamma c-\beta b)>0.
\]
It follows that \(I(\alpha,\beta,\gamma)\subset(0,\infty)\) for all positive parameters.

Next consider the edge at \(\beta b\), and define
\[
K(\alpha,\beta,\gamma)
=\{\varphi(\gamma c-\beta b):\varphi\in J(\beta b-\alpha a)\}.
\]
This interval also avoids zero. Indeed, if \(\varphi(\gamma c-\beta b)=0\), then \(-\varphi\in J(\alpha a-\beta b)\) and \(\alpha a-\beta b\perp_B\gamma c-\beta b\). The restricted property for the same rescaled triple would then be violated. Its sign is therefore constant on the positive octant. To determine that sign, fix \(\beta>0\) and let \(\alpha,\gamma\) be sufficiently small. Every cluster point of \(J(\beta b-\alpha a)\) as \(\alpha\downarrow0\) belongs to \(J(b)\), and each such functional takes the strictly positive value \(N(b)\) on \(b\). Compactness therefore gives a uniform positive lower bound for \(\varphi(b)\) when \(\alpha\) is small. Taking \(\gamma\) still smaller makes
\[
\varphi(\gamma c-\beta b)=\gamma\varphi(c)-\beta\varphi(b)<0
\]
for every \(\varphi\in J(\beta b-\alpha a)\). Thus \(K\subset(-\infty,0)\) throughout the positive octant.

Fix \(\varphi\in J(\beta b-\alpha a)\). The sign conclusions give, for every \(\gamma>0\),
\[
\varphi(\gamma c-\alpha a)>0,\qquad
\varphi(\gamma c-\beta b)<0.
\]
Letting \(\gamma\to\infty\) yields \(\varphi(c)\geq0\) and \(\varphi(c)\leq0\), so \(\varphi(c)=0\). Substitution into the strict inequalities gives \(\varphi(a)<0<\varphi(b)\). Relabeling \(a,b,c\) proves \eqref{eq:annihilation}-\eqref{eq:strict-signs}.
\end{proof}

No two residuals in a triple satisfying \eqref{eq:S} are collinear. If two were negative multiples, the inequalities \(p(u_i)<0\) would fail for one of them. If \(u_j\) and \(u_k\) were positive multiples, choose balancing functionals as in \eqref{eq:ft-balance} and evaluate them at their common unit direction. The first two functionals would each take the value \(1\), forcing the third to take the value \(-2\), which is impossible for a unit functional.

\section{Nonsmooth points and their geometric properties}
\label{sec:nonsmooth}

\begin{lemma}\label{lem:nonsmooth}
Every residual of a triple satisfying \eqref{eq:S} is a nonsmooth point of \(N\).
\end{lemma}

\begin{proof}
By Lemma~\ref{lem:radial}, normalize each residual to have norm one. Choose \(f\in J(u_1)\), \(g\in J(u_2)\), and \(h\in J(u_3)\) with \(f+g+h=0\). Suppose \(u_1\) is smooth, so \(J(u_1)=\{f\}\). Apply \eqref{eq:annihilation} to \(u_1-tu_2\) and \(u_1-tu_3\), and let \(t\downarrow0\). Closedness of the graph of \(J\) gives
\[
f(u_2)=f(u_3)=0.
\]
The balance relation now yields \(g(u_3)=h(u_2)=-1\). For every \(t>0\),
\[
g(u_2-tu_3)=1+t=N(u_2-tu_3),
\]
where the equality follows from \(g\in J(u_2)\), the value \(g(u_3)=-1\), and the triangle inequality. Thus \(g\in J(u_2-tu_3)\), and \eqref{eq:annihilation} gives \(g(u_1)=0\). The same argument gives \(h(u_1)=0\). But then
\[
(f+g+h)(u_1)=f(u_1)=1,
\]
contradicting \(f+g+h=0\).
\end{proof}

Let
\[
\Sigma=\{x\in\R^3\setminus\{0\}:J(x)\text{ is not a singleton}\}.
\]

\begin{lemma}\label{lem:graph-cover}
There are countably many global Lipschitz maps \(R_\ell:\R^2\to\R^3\) whose images cover \(\Sigma\). At every differentiability point, \(DR_\ell\) has rank two.
\end{lemma}

\begin{proof}
Use Euclidean coordinates also on the dual. Choose rational \(a,b\in\R^3\) and rational \(\delta>0\) with \(2\delta<|a-b|\). Let \(S(a,b,\delta)\) consist of the nonzero \(x\) for which \(J(x)\) contains \(f,g\) with \(|f-a|<\delta\) and \(|g-b|<\delta\).

For \(x,y\in S(a,b,\delta)\), choose such \(f_x,g_x,f_y,g_y\). Monotonicity of the subdifferential of \(N\) gives
\[
(f_x-g_y)(x-y)\geq0,\qquad (g_x-f_y)(x-y)\geq0,
\]
and hence
\[
|(a-b)(x-y)|\leq2\delta|x-y|.
\tag{3.1}\label{eq:cone-estimate}
\]
Set \(e=(a-b)^\sharp/|a-b|\), where \((a-b)^\sharp\) is the Euclidean vector representing the covector \(a-b\), and put \(\kappa=2\delta/|a-b|<1\). Estimate \eqref{eq:cone-estimate} shows that orthogonal projection onto \(e^\perp\) is injective on \(S(a,b,\delta)\). It also shows that the \(e\)-coordinate is a Lipschitz function of the projected coordinate, with constant at most \(\kappa/\sqrt{1-\kappa^2}\). Extend this scalar function from its domain in \(e^\perp\) to all of \(e^\perp\) by the McShane extension theorem. Its graph defines a global Lipschitz map \(R:\R^2\to\R^3\), and its horizontal projection is the identity, so its differential has rank two wherever it exists.

Every nonsmooth point admits two distinct norming functionals. Rational centers and a sufficiently small rational \(\delta\) place these two functionals in one of the sets \(S(a,b,\delta)\). The family of rational choices is countable, proving the cover.
\end{proof}

By Rademacher's theorem \cite[Theorem~2.5.4]{BV2010}, \(R_\ell\) and \(N\circ R_\ell\) are differentiable almost everywhere. At such a parameter \(t\), put \(u=R_\ell(t)\) and \(H=\operatorname{ran}DR_\ell(t)\). For any \(f\in J(u)\), the support inequality applied to \(R_\ell(t\pm sh)\) implies
\[
D(N\circ R_\ell)(t)h=f(DR_\ell(t)h)\qquad(h\in\R^2).
\]
Consequently,
\[
J(u)-J(u)\subset H^\perp.
\tag{3.2}\label{eq:normal-differences}
\]
If \(u\in\Sigma\), the left side contains a nonzero functional. Since \(H\) is a plane, the linear span of \(J(u)-J(u)\) is exactly its one-dimensional annihilator.

We next establish an estimate near points where Alexandrov's theorem gives a second-order expansion.

\begin{lemma}\label{lem:alex-estimate}
There is a Borel Lebesgue-null set \(B\subset\R^3\), containing the origin, such that for each \(q\notin B\), \(N\) has a second-order expansion
\[
N(q+h)=N(q)+\ell_q(h)+\tfrac12 A_q[h,h]+o(|h|^2),
\tag{3.3}\label{eq:alex-expansion}
\]
and there are \(C_q,\rho_q>0\) for which
\[
|\eta-\ell_q|\leq C_q|y-q|
\quad\bigl(0<|y-q|<\rho_q,\ \eta\in J(y)\bigr).
\tag{3.4}\label{eq:all-subgradients}
\]
In particular, \(J(q)=\{\ell_q\}\).
\end{lemma}

\begin{proof}
By Alexandrov's theorem \cite[Theorems~2.6.1 and 2.6.4]{BV2010}, the expansion \eqref{eq:alex-expansion} holds outside a null set. It implies differentiability at \(q\), so \(J(q)=\{\ell_q\}\). Put \(h=y-q\) and \(r=|h|\). For each vector \(e\) in a Euclidean unit basis, the supporting inequality at \(y\), applied to \(y+re\) and \(y-re\), gives
\[
\frac{N(y)-N(y-re)}r\leq\eta(e)\leq
\frac{N(y+re)-N(y)}r.
\]
Substituting \eqref{eq:alex-expansion} into both bounds yields
\[
|\eta(e)-\ell_q(e)|\leq C r
\]
for each basis vector, with a constant independent of \(\eta\). Since there are only three basis vectors, the Euclidean dual norm satisfies \eqref{eq:all-subgradients}.
\end{proof}

\section{Proof of the main theorem}
\label{sec:selection}

\begin{proof}[Proof of Theorem~\ref{thm:main}]

Suppose that \(N\) is not induced by an inner product, while the restricted property \eqref{eq:restricted-property} holds. By \cite[Theorem~15]{BFS2002}, there is a triple \(A^0=(A_1^0,A_2^0,A_3^0)\) such that
\[
Z_1(A_1^0,A_2^0,A_3^0)\cap\co\{A_1^0,A_2^0,A_3^0\}=\varnothing.
\]
In finite dimensions the median set is nonempty, compact, and convex. Strict separation therefore gives a nonzero linear functional \(p\) such that
\[
\min_{m\in Z_1(A_1^0,A_2^0,A_3^0)}p(m)>
\max_i p(A_i^0).
\tag{4.1}\label{eq:strict-separation}
\]

The median map is upper semicontinuous and locally bounded. To see this, note that
\[
\sum_{i=1}^3N(A_i-z)\geq3N(z)-\sum_{i=1}^3N(A_i)
\]
gives a common bound on nearby minimizers, and continuity of the objective then gives upper semicontinuity. After shrinking to a bounded open neighborhood \(U\subset(\R^3)^3\) of \(A^0\), we therefore have
\[
\min_{m\in Z_1(A_1,A_2,A_3)}p(m)>
\max_i p(A_i)\qquad(A\in U).
\tag{4.2}\label{eq:open-failure}
\]
For each \(A\in U\), choose a median \(m\) at which \(p\) attains its minimum over the median set. Then the residuals \(u_i=A_i-m\) satisfy \eqref{eq:S}. By Lemma~\ref{lem:nonsmooth}, each \(u_i\) belongs to \(\Sigma\).

Choose three maps \(R_1,R_2,R_3\) from the countable cover in Lemma~\ref{lem:graph-cover} and define
\[
F(m,t_1,t_2,t_3)
=\bigl(m+R_1(t_1),m+R_2(t_2),m+R_3(t_3)\bigr),
\qquad
F:\R^9\to\R^9.
\tag{4.3}\label{eq:F}
\]
For each such chart triple, let \(E\) be the set of parameters for which \(F(m,t_1,t_2,t_3)\in U\) and \(u_i=R_i(t_i)\) satisfy \eqref{eq:S}.

The set \(E\) is Borel. For the verification, define
\[
D_u(z)=\sum_iN(u_i-z),\qquad V(u)=\min_zD_u(z),\qquad
V_n(u)=\min_{p(z)\leq-1/n}D_u(z).
\]
Coercivity ensures that all these minima are attained. The functions \(V\) and \(V_n\) are Lipschitz in \(u\): changing \(u\) to \(u'\) changes every value of the objective by at most \(\sum_iN(u_i-u_i')\). Condition \eqref{eq:S} is equivalent to
\[
D_u(0)=V(u),\qquad V_n(u)>D_u(0)\ (n\geq1),\qquad p(u_i)<0\ (i=1,2,3).
\tag{4.4}\label{eq:Borel-S}
\]
These are Borel conditions. As the chart triple varies over a countable family, the sets \(F(E)\) cover \(U\).

We use the equal-dimensional area formula for Lipschitz maps:
\[
\int_A|\det DG(x)|\,dx
=\int_{\R^9}\#\bigl(A\cap G^{-1}(y)\bigr)\,dy
\tag{4.5}\label{eq:area-formula}
\]
for measurable \(A\subset\R^9\) and Lipschitz \(G:\R^9\to\R^9\); see \cite[Sec.~3.3]{EG2015}. For each chart triple, discard the null set where one of the maps \(R_i\) or \(N\circ R_i\) fails to be differentiable. Its image under \(F\) is null. The area formula also shows that the image of the critical set \(\{\det DF=0\}\) has measure zero. The sets \(F(E)\), over the countable family of chart triples, cover the open set \(U\). Hence for at least one chart triple, the following set has positive measure:
\[
\begin{aligned}
E_+=\{x\in E:\;&\det DF(x)\ne0,\quad
R_i\text{ and }N\circ R_i\text{ are differentiable at }t_i\\
&\text{for }i=1,2,3\}.
\end{aligned}
\tag{4.6}\label{eq:Eplus}
\]
Fix a point of \(E_+\). Write \(H_i=\operatorname{ran}DR_i(t_i)\), and choose a nonzero normal covector \(\psi_i\) to each plane. The derivative of \(F\) is invertible precisely when \(\psi_1,\psi_2,\psi_3\) are linearly independent. Indeed, a vector in the kernel must have its \(m\)-component in \(H_1\cap H_2\cap H_3\); conversely, a common vector in these planes determines a kernel vector by solving in each two-dimensional chart.

For positive rational numbers \(a_1,a_2,a_3\), define
\[
F_a(m,t_1,t_2,t_3)
=\bigl(m+a_1R_1(t_1),m+a_2R_2(t_2),m+a_3R_3(t_3)\bigr).
\]
At differentiability points,
\[
|\det DF_a|=a_1^2a_2^2a_3^2|\det DF|.
\tag{4.7}\label{eq:scaled-jacobian}
\]
For each pair \(j\ne k\), the set of target triples with \(A_j-A_k\in B\) has measure zero in \(\R^9\), by a linear change of variables and Fubini's theorem. Applying the area formula to \(F_a\) shows that the corresponding preimage in \(E_+\) is null, since the Jacobian is positive there. Only countably many positive rational choices of \(a_i\) and pairs \(j,k\) occur. We can therefore choose \(x_0=(m,t_1,t_2,t_3)\in E_+\) that is also a density point of \(E\) and satisfies
\[
a_j u_j-a_k u_k\notin B
\quad\text{for every }j\ne k\text{ and every }a_j,a_k\in\Q_{>0},
\qquad u_i=R_i(t_i).
\tag{4.8}\label{eq:good-mixed}
\]

Fix distinct \(i,j,k\) and positive rational \(s,t\). Set \(q=su_j-tu_k\). By \eqref{eq:good-mixed} and Lemma~\ref{lem:alex-estimate}, \(J(q)=\{\varphi\}\). Let \(h\in\R^2\) be arbitrary and perturb only the \(i\)-th chart parameter. Since \(x_0\) is a density point of \(E\), there are \(\varepsilon_n\downarrow0\) and \(x_n\in E\) such that
\[
x_n=x_0+\varepsilon_n(0,\ldots,h,\ldots,0)+o(\varepsilon_n).
\tag{4.9}\label{eq:density-perturbation}
\]
To obtain this sequence, center balls of radius \(o(\varepsilon_n)\) at the indicated points. Since \(x_0\) is a density point of \(E\), the complement of \(E\) has vanishing relative measure in balls of radius \(O(\varepsilon_n)\). The smaller balls can therefore be chosen to intersect \(E\).

Write \(u_{\ell,n}=R_\ell(t_{\ell,n})\). Differentiability of the \(i\)-th chart and Lipschitz continuity of the others give
\[
u_{i,n}=u_i+\varepsilon_nDR_i(t_i)h+o(\varepsilon_n),\qquad
u_{j,n}=u_j+o(\varepsilon_n),\qquad
u_{k,n}=u_k+o(\varepsilon_n).
\tag{4.10}\label{eq:residual-variation}
\]
Hence \(q_n=su_{j,n}-tu_{k,n}=q+o(\varepsilon_n)\). Choose any \(\varphi_n\in J(q_n)\). Estimate \eqref{eq:all-subgradients} at \(q\) yields \(\varphi_n=\varphi+o(\varepsilon_n)\). Since the residuals at \(x_n\) satisfy \eqref{eq:S}, Lemma~\ref{lem:mixed-edge} gives
\[
\varphi_n(u_{i,n})=0,\qquad \varphi(u_i)=0.
\]
Subtracting, dividing by \(\varepsilon_n\), and using \eqref{eq:residual-variation} shows that \(\varphi(DR_i(t_i)h)=0\). Since \(h\) is arbitrary,
\[
J(su_j-tu_k)\subset H_i^\perp=\operatorname{span}\{\psi_i\}
\qquad(s,t\in\Q_{>0}).
\tag{4.11}\label{eq:edge-normal}
\]

For each of the cyclic choices \((i,j,k)=(1,2,3),(2,3,1),(3,1,2)\), denote the unique element of \(J(u_j-u_k)\) by \(\varphi_i\). Every element of \(J(su_j-tu_k)\) has dual norm one and, by \eqref{eq:edge-normal}, lies on the line spanned by \(\varphi_i\); it is therefore either \(\varphi_i\) or \(-\varphi_i\). The strict sign \eqref{eq:strict-signs} fixes the positive choice, so
\[
J(su_j-tu_k)=\{\varphi_i\}\qquad(s,t\in\Q_{>0}).
\tag{4.12}\label{eq:all-edge-supports}
\]
The three functionals \(\varphi_i\) are linearly independent because their respective lines are the three independent normal directions \(\operatorname{span}\{\psi_i\}\).

Letting the rational ratios in \eqref{eq:all-edge-supports} tend to zero or infinity and using the closed graph of \(J\), we obtain
\[
\begin{array}{ll}
\varphi_1\in J(u_2),&-\varphi_1\in J(u_3),\\
\varphi_2\in J(u_3),&-\varphi_2\in J(u_1),\\
\varphi_3\in J(u_1),&-\varphi_3\in J(u_2).
\end{array}
\tag{4.13}\label{eq:cyclic-supports}
\]
However, \eqref{eq:normal-differences} applied at \(u_1\) gives
\[
J(u_1)-J(u_1)\subset\operatorname{span}\{\psi_1\}
=\operatorname{span}\{\varphi_1\}.
\]
Since \(\varphi_3,-\varphi_2\in J(u_1)\), it follows that
\[
\varphi_2+\varphi_3\in\operatorname{span}\{\varphi_1\},
\]
contradicting the linear independence of \(\varphi_1,\varphi_2,\varphi_3\). Thus every three-dimensional norm satisfying the restricted property is induced by an inner product.

We now return to arbitrary dimension. Let \(Y\subset X\) be any three-dimensional subspace. The restricted property passes from \(X\) to \(Y\): for \(x,y\in Y\), the point supplied by the property lies in \(\co\{0,x,y\}\subset Y\) and is a median for the norm restricted to \(Y\). The three-dimensional result shows that every such subspace is an inner-product space. Since any pair of vectors in \(X\) lies in a three-dimensional subspace, the parallelogram identity holds for every pair in \(X\). The polarization identity then defines an inner product that induces the norm.

For the converse, suppose the norm comes from an inner product. Fix any triple \(a_1,a_2,a_3\), and put \(C=\co\{a_1,a_2,a_3\}\). If \(z\notin C\), let \(c\) be the metric projection of \(z\) onto \(C\). The projection inequality gives \(\langle a_i-c,z-c\rangle\leq0\) for each \(i\), and hence
\[
\|a_i-z\|^2
=\|a_i-c\|^2+\|z-c\|^2-2\langle a_i-c,z-c\rangle
>\|a_i-c\|^2.
\]
Thus \(z\) cannot minimize the sum of distances. Every Fermat-Torricelli median lies in \(C\), and in particular the restricted property holds. Theorem~\ref{thm:main} follows.
\end{proof}

\section*{Acknowledgements}

The authors thank the members of the Functional Analysis Research team at the School of Mathematics and Statistics of Anqing Normal University for their helpful discussions and for pointing out errors in earlier versions of this paper. 

\section*{Declarations}

\subsection*{Competing Interests}
The authors declare that there are no competing interests.

\subsection*{Funding Information}
Not applicable.

\subsection*{Author contribution}
The authors contributed equally to this work. Both authors have read and approved the manuscript.

\subsection*{Data Availability Statement}
Not applicable.

\end{document}